\documentclass[11pt]{amsart}

\usepackage{amsmath,amsthm, amscd, amssymb, amsfonts}
\usepackage[all]{xy}
\usepackage[dvips, dvipsnames, usenames]{color}
\allowdisplaybreaks
\usepackage{enumitem}

\usepackage{xcolor}

\newtheorem{theorem}{Theorem}[section]

\newtheorem{lemma}[theorem]{Lemma}

\newtheorem{teo intro}{Theorem}

\newtheorem{proposition}[theorem]{Proposition}

\theoremstyle{definition}
\newtheorem{definition}[theorem]{Definition}

\theoremstyle{remark}
\newtheorem{remark}[theorem]{Remark}

\newcommand\pf{\begin{proof}}
	\newcommand\epf{\end{proof}}

\newcommand\yd{{}^{H}_{H}\mathcal{YD}}

\newcommand\Rat{\operatorname{Rat}}

\newcommand\ot{\otimes}
\newcommand\ra{\rightarrow}

\newcommand\N{\mathbb{N}}

\newcommand\cG{\mathcal{G}}
\newcommand{\cF}{\mathcal{F}}

\newcommand\cS{\mathcal{S}}

\newcommand{\ku}{ \mathbf{k}}

\newcommand{\M}{\mathcal M}

\begin{document}
\title[On braided co-Frobenius Hopf algebras]{On braided co-Frobenius Hopf algebras}
\author[Rossi Bertone]
{Fiorela Rossi Bertone}

\address{Instituto de Matem\'atica (INMABB), Departamento de Matem\'atica, Universidad Nacional del Sur (UNS)-CONICET, Bah\'ia Blanca, Argentina
} \email{fiorela.rossi@uns.edu.ar}

\thanks{\noindent 2000 \emph{Mathematics Subject Classification.}
	16W30.}

\dedicatory{To Nicol\'as Andruskiewitsch on his 60th birthday.}

\begin{abstract}
	We present characterizations of braided co-Frobenius Hopf algebras in the braided tensor category of Yetter-Drinfeld modules over a Hopf algebra extending those already known for co-Frobenius Hopf algebras.
\end{abstract}

\maketitle

\section*{Introduction}
The Haar measure on a compact Lie group can be spelled out in terms of the algebra of rational functions on the group as shown by Hochschild. Sweedler extended this definition to arbitrary Hopf algebras \cite{Sw}.
Namely, a \emph{left integral} on a Hopf algebra $H$ is $\int\in H^*$ satisfying $\alpha\cdot \int=\alpha(1_H)\int$ for all $\alpha\in H^*$, or equivalently, $\int(h_{(2)})h_{(1)}=\int(h)1_H$ for all $h\in H$
(analogously, one  defines right integrals). The linear subspace of $H^*$ formed by all left integrals is denoted $\int_{\ell}(H)$. Then $\dim \int_{\ell} (H) \leq 1$ \cite[p. 331]{Sw}. If $\int_{\ell}(H) \neq 0$, then $H$ is \emph{co-Frobenius}.
Co-Frobenius Hopf algebras have been intensively studied and several equivalent characterizations have been presented in the literature, see Theorem \ref{thm:characterizations_coFrob} for details and references. 
Some of these characterizations have a cohomological flavour. Another states that co-Frobenius is equivalent to the 
finiteness of the coradical filtration \cite{R-cofrob,AC,ACE}; see the Introductions of these articles for detailed discussions.

The purpose of this paper is to extend several of these characterizations to Hopf algebras $R$ in the braided tensor category $\yd$ of Yetter-Drinfeld modules over $H$. Our main results are Theorems \ref{Thm. cofrob eq for R} and 
\ref{Thm. cofrob eq for R-Hcofrob} (the second requires that $H$ is co-Frobenius as well). 
As for the relation with the finiteness of the coradical filtration of $R$, one implication does not present difficulties, see
Lemma \ref{lem:finfilt_cofrob}. The converse is only proved in specific conditions and remains an open question in general. See 
Proposition \ref{prop:cofrob-finite-corad-filt}.

\subsection*{Acknowledgements}
This work is part of author's PhD Thesis. The author wants to thank her advisor Nicol\'as Andruskiewitsch for his guidance during these years, the careful reading and his many suggestions on this work which significantly improved it.

\section{Preliminaries}
\subsection{Notations}

Let $\ku$ be a field. All vector spaces, tensor products, algebras, coalgebras, 
are over $\ku$.
We use the Sweedler notation for coalgebras and comodules; the counit and the antipode of a Hopf algebra (always assumed bijective in this paper) are denoted by $\varepsilon$ and $\cS$ respectively. 
Our reference for Hopf algebras is \cite{Rad-libro}.
The coradical filtration of a coalgebra $C$ is denoted by $(C_n)_{n\in \N_0}$.

Throughout the article $H$ is a Hopf  algebra.
An $H$-\emph{module algebra} is an algebra in the monoidal category $_H\hspace{-2pt}\M$ of $A$-modules, while  a (left) $H$-\emph{comodule coalgebra}  is a coalgebra in the category $^H\hspace{-2pt} \M$ of (left) $H$-comodules.
If $R$ is a left $H$-module algebra, then the \emph{smash product} of $R$ and $H$ is the vector space $R\otimes H$ with the multiplication 
	\begin{align}
	(r\ot h)(s\ot k)=r(h_{(1)}.s)\ot h_{(2)}k, \qquad r,s\in R, \, h,k\in H.
	\end{align}
As well, if $R$ is a left $H$-comodule coalgebra, the \emph{smash coproduct} of $R$ and $H$ is the vector space $R\otimes H$ with the comultiplication 
	\begin{align}
	\Delta(s\ot h)=s^{(1)}\ot (s^{(2)})^{(-1)}h_{(1)} \otimes (s^{(2)})^{(0)}\ot h_{(2)}.	
	\end{align}
By abuse of notation,  $R\#H$ denotes either the smash product or coproduct.
A \emph{Yetter-Drinfeld module} over $H$ is a vector space $V$ which is an $H$-module and an $H$-comodule satisfying the compatibility condition
\begin{align*}
\rho(h\cdot v)=h_{(1)}v_{(-1)}\cS(h_{(3)})\otimes h_{(2)}\cdot v_{(0)}, & &\forall h\in H, v\in V.
\end{align*}
The category $\yd$ of Yetter-Drinfeld modules over $H$ is a (braided) tensor category. Thus there are Hopf algebras in $\yd$
(also called braided Hopf algebras); for them, we emphasize the braided structure using the notation $\underline{\Delta}$ for the coproduct and a variation of the Sweedler notation $\underline{\Delta}(x)=x^{(1)}\ot x^{(2)}$. 
The bosonization of $R$ by $H$  is the  Hopf algebra $R\# H$ with the smash product and smash coproduct \cite{Ma,R}.

\subsection{Equivalent monoidal categories}

Let $R$ be an algebra in $_H \M$ 
and $S$ a coalgebra in $^H\hspace{-2pt} \M$.
We denote by $_R\underline{\M}$ the category of left $R$-modules inside $_H \M$ (i.e.\ for $V\in {}_R\underline{\M}$, the $R$-action on $V$ is a morphism of $H$-modules) and  by $^S\hspace{-2pt}\underline{\M}$  
the category of left $S$-comodules inside $^H\hspace{-2pt} \M$. 
The next proposition is wellknown.

\begin{lemma}\label{prop: equivalent_cat}
There are equivalences of abelian categories
\begin{align}\label{eq:equivalent_cat}
_{R\# H}\M &\simeq {}_R\underline{\M}, &
^{S\#H}\hspace{-2pt} \M&\simeq {}^S\hspace{-2pt}\underline{\M}.
\end{align}
If $R$, respectively $S$, is a Hopf algebra in $\yd$, then the equivalences in \eqref{eq:equivalent_cat} are actually equivalences of monoidal categories.
\end{lemma}
\pf
Let $V\in {}_R\underline{\M}$ with actions $\triangleright: H \otimes V\to V$,
$\rightharpoonup: R\otimes V\to V$ such that  
\begin{align}\label{eq: compat Rmod}
h\triangleright(r\rightharpoonup v)=(h_{(1)}.r)\rightharpoonup (h_{(2)}\triangleright v), \qquad \forall h\in H,\, r\in R,\, v\in V.
\end{align}
Let us define $\cdot:R\# H\otimes V \longrightarrow V$ by $(r\# h)\cdot v := r\rightharpoonup(h\triangleright v)$, then $V \in {}_{R\# H}\M $ by  \eqref{eq: compat Rmod}. 
 Thus, we have the functor $\cG:{}_R\underline{\M}\longrightarrow{}_{R\# H} \M$ which is the identity on morphisms.
Conversely, the  morphisms $\iota: H \rightarrow R\# H$, $h\mapsto 1\# h$ and $\pi: R\# H\rightarrow H$, $r\#h\mapsto \varepsilon(r)h$ give rise to the funtor $\cF:{}_A\M\longrightarrow {}_R\underline{\M}$, $\cF(V)=V$ with the module structure given by
$$r\rightharpoonup v:= (r\# 1)\cdot v  \qquad \mbox{and} \qquad
h\triangleright v:= (1\# h)\cdot v.
$$
It is easy to see that $\cF$ and $\cG$ are mutually inverse. 
The first statement follows.
The proof of the second is analogous, while the claims on monoidal equivalences are straightforward. See \cite{RB-tesis} for details.
\epf

\subsection{Co-Frobenius Hopf algebras}
Here we collect the  basic definitions and various characterizations of co-Frobenius Hopf algebras.

Let $C$ be a coalgebra. The maximal rational submodule of the left (respectively right) $C^*$-module $C^*$ is denoted $\Rat(C^*)$ (respectively $\Rat(C^*_{C^*})$).
Also, $E_C (M)$ stands for the injective hull of $M \in {}^C\M$ (or 
$\M^C$).

\begin{definition}\cite{DNR}
The coalgebra $C$  is (left)
\begin{itemize}[leftmargin=*]
	\item  \emph{semiperfect}  if  $\dim E_C(S) < \infty$ for any   $S\in {}^C\M$ simple;
	\item \emph{quasi co-Frobenius} if there exist a free $C^*$-module $F$ and an injective morphism of left $C^*$-modules $f:C\ra F$;
	\item \emph{co-Frobenius} if there exists a monomorphism 
	$C \to C^*$ of left $C^*$-modules.
\end{itemize}
\end{definition}

Clearly, a Hopf algebra $A$ is (left) \emph{co-Frobenius} iff $\int_{\ell}(A) \neq 0$. Next, we list some equivalent formulations for a Hopf algebra $A$.

\begin{theorem}\label{thm:characterizations_coFrob}
The following assertions are equivalent:
\begin{enumerate}[leftmargin=*,label=\rm{($\roman*$)}]
\item $A$ is co-Frobenius.
\item $E_A(S)$ is finite-dimensional for all simple $S\in {}^A\hspace{-2pt}\M$.
\item $E_A(\ku)$ is finite-dimensional.
\item There exists a non-zero injective $A$-comodule with finite dimension.
\item  Every $0 \neq M \in {}^A\hspace{-2pt}\M$ has a finite-dimensional non-zero quotient.
\item Every $0 \neq M \in {}^A\hspace{-2pt}\M$ injective has a finite-dimensional $\neq 0$ quotient.
\item There is $0 \neq M \in {}^A\hspace{-2pt}\M$ injective with a finite-dimensional $\neq 0$ quotient.
\item ${}^A\hspace{-2pt}\M$ has non-zero projective objects.
\item Every  $M\in ^A\hspace{-2pt}\M$ has a projective cover.
\item Every injective right $A$-comodule is projective.
\item $\Rat(A^*) \neq 0$.
\item The coradical filtration of $A$ is finite.
\end{enumerate}
\end{theorem}
\pf
$(i)\Leftrightarrow (ii)$ is \cite[Theorem 3]{L} and $(ii)\Leftrightarrow(ix) \Leftrightarrow (xi)$ is \cite[Theorem 10]{L}. $(ii)\Rightarrow(iii)\Rightarrow(iv)$ is trivial, see \cite[Footnote 1]{AC}. $(iv)\Rightarrow(i)$ is \cite[Theorem 2.3]{DN}. $(i)\Leftrightarrow(v)\Leftrightarrow(vi)\Leftrightarrow(vii)$ is \cite[Theorem 2.3]{AC} and $(i)\Leftrightarrow(viii)\Leftrightarrow(x)$ is \cite[Theorem 2.8]{AC}. $(xii) \implies (i)$ is \cite[Theorem 2.1]{AD} and $(i) \implies (xii)$ is  \cite[Theorem 2]{ACE}.
\epf

\section{Braided co-Frobenius Hopf algebras}
\subsection{First characterizations}
Let $R$ be a Hopf algebra  in $\yd$ and  $A=R\# H$.
A right  integral on $R$ is a map $\lambda: R\to \ku$ such that $\lambda f=f(1)\lambda$  for all $f\in R^*$.
A braided Hopf algebra $R\in \yd$ is called right  co-Frobenius if it admits a non-zero right  integral; 
equivalently, the coalgebra underlying $R$ is right co-Frobenius.
The space of right integrals on $R$ is denoted $\int_r(R)$.
Left versions are defined and denoted in the obvious way.
We derive now generalizations of Theorem \ref{thm:characterizations_coFrob} for Hopf algebras in $\yd$.

Since a braided version of the fundamental theorem of Hopf modules holds in this context \cite[Thm. 3.4]{Takeuchi}, \cite[Lemma 2.1]{MS} and $\Rat(R^*_{R^*})$ is a (braided) Hopf module \cite[\S 2]{AG}, we have the following   result, see also \cite[\S 2]{GZ}.

\begin{lemma}\label{thm:rat_int}
 $\int_r(R)\ot R\simeq \Rat(R^*_{R^*})$;
hence $\Rat(R^*_{R^*})\neq0$ iff $\int_r (R) \neq0$.\qed
\end{lemma}

\begin{remark}\label{rmk:left-right coFrob}
If $R$ is left co-Frobenius then it is right co-Frobenius.
Indeed, $R$ is a left co-Frobenius coalgebra, hence  a left quasi co-Frobenius coalgebra. Thus, $R$ is left semiperfect  and $\Rat(R^*_{R^*})\neq0$, see e.~g. \cite[Corollary 3.3.6]{DNR}. By Lemma \ref{thm:rat_int}, $R$ has a non-zero right integral and the claim follows.
\end{remark}

Here is a version of Theorem \ref{thm:characterizations_coFrob} for $R\in \yd$. 

\begin{theorem}\label{Thm. cofrob eq for R}
The following assertions are equivalent:
\begin{enumerate}[leftmargin=*,label=\rm{($\roman*$)}]
\item $R$ is co-Frobenius.
\item $E_R(S)$ is finite-dimensional for all simple  $S\in{}^R\hspace{-2pt}\M$.
\item  Every  $0 \neq M \in {}^R\hspace{-2pt}\M$ has a finite-dimensional non-zero quotient.
\item Every $0 \neq M \in {}^R\hspace{-2pt}\M$ injective  has a finite-dimensional $\neq 0$ quotient.
\item There is $0 \neq M \in {}^R\hspace{-2pt}\M$ injective with a finite-dimensional $\neq 0$ quotient. 
\item Every object in ${}^R\hspace{-2pt}\M$ has a projective cover.
\item Every injective object in ${}^R\hspace{-2pt}\M$ is projective.
\item $\Rat(R^*) \neq 0$.
\end{enumerate}
\end{theorem}
\pf
$(i) \Rightarrow (ii)$: co-Frobenius implies quasi co-Frobenius  which implies semiperfect \cite[3.3.6]{DNR}.
By \cite[3.2.3]{DNR},   $(ii) \Leftrightarrow (vi)$. By \cite[3.3.4, 3.3.6]{DNR} $ (i) \Rightarrow (vii)\Rightarrow (viii)$; Lemma \ref{thm:rat_int} gives $(i)\Leftrightarrow (viii)$.  
The implications $(iii)\Rightarrow(iv)\Rightarrow(v)$ are obvious. 
Finally, $(ii)\Rightarrow(iii)$ and  $(v)\Rightarrow(viii)$ by the same proofs of $(i)\Rightarrow(ii)$ and  $(iv)\Rightarrow(i)$ in \cite[Theorem 2.3]{AC}.
\epf

\begin{lemma}\label{lem:finfilt_cofrob}
If the coradical filtration of $R$ is finite, then $R$ is co-Frobenius.
\end{lemma}

\pf
By \cite[Theorem 1.2]{AD}, $\Rat(R^*)\neq0$. Then, by the left version of Lemma \ref{thm:rat_int}, $R$ is co-Frobenius.
\epf

\subsection{When $H$ is co-Frobenius}
\begin{lemma}\label{prop:RcoF_Acof}
The Hopf algebra $A=R\# H$
is co-Frobenius if and only if both $H$ and $R$ are so.
\end{lemma}
\pf
If $H$ and $R$ are co-Frobenius, pick $0 \neq \lambda_R\in \int_r(R)$ and $0 \neq \lambda_H \in \int_r(H)$.
By \cite[Proposition 4]{R} $\lambda_R\#\lambda_H \in \int_r(R \# H)$ and $A$ is co-Frobenius. 

Assume that $A$ is co-Frobenius. Then so is its Hopf subalgebra $H$
\cite[Theorem 2.15]{Su}. Let $0 \neq \lambda \in \int_{\ell}(A)$. 
Fix $r\#h\in A$ such that $\lambda(r\#h)\neq0$ and set $\lambda_R(s)=\lambda(s\# h)$, $s\in R$. Clearly $\lambda_R \neq 0$.
We claim that $\lambda_R \in \int_{\ell}(R)$.  Let $\psi:A\rightarrow R$, $\psi(s\#k)=\varepsilon(k)s$. Then, for all $s\in R$, we have
\begin{align*}
\lambda_R(s)1_R&=\psi\circ\iota(\lambda_R(s)1_R)=
\psi(\lambda(s\#h)1_A)=\psi((s\#h)_{(1)}\lambda((s\#h)_{(2)})) \\
&=\psi(s_{(1)}\#(s_{(2)})^{(-1)}h_{(1)})\lambda((s_{(2)})^{(0)}\#h_{(2)})\\
&=
\varepsilon((s_{(2)})^{(-1)}h_{(1)})\, s_{(1)}\lambda((s_{(2)})^{(0)}\#h_{(2)})
=s_{(1)}\lambda_R(s_{(2)}), 
\end{align*}
and the theorem follows.
\epf

Theorem  \ref{Thm. cofrob eq for R} admits a refinement.

\begin{theorem}\label{Thm. cofrob eq for R-Hcofrob}
If $H$ is a co-Frobenius Hopf algebra,  then the  following assertions are equivalent for a braided Hopf algebra $R\in\yd$:
\begin{enumerate}[leftmargin=*,label=\rm{($\roman*$)}]
\item $R$ is co-Frobenius.
\item $E_R(S)$ is finite-dimensional for all simple object in $S\in{}^R\hspace{-2pt}\underline{\M}$.
\item $E_R(\ku)$ is finite-dimensional.
\item There exists a non-zero injective object in $^R\hspace{-2pt}\underline{\M}$ with finite dimension.
\item  Every  $0 \neq M \in {}^R\hspace{-2pt}\underline{\M}$ has a finite-dimensional $\neq 0$ quotient.
\item Every $0 \neq M \in {}^R\hspace{-2pt}\underline{\M}$ injective has a finite-dimensional $\neq 0$ quotient.
\item There is $0 \neq M \in {}^R\hspace{-2pt}\underline{\M}$ injective with a finite-dimensional $\neq 0$ quotient.
\item $^R\hspace{-2pt}\underline{\M}$ has non-zero projective objects.
\item Every comodule $V\in{}^R\hspace{-2pt}\underline{\M}$ has a projective cover.
\item Every injective object in $^R\hspace{-2pt}\underline{\M}$ is projective.
\end{enumerate}
\end{theorem}
\pf
By Lemma \ref{prop:RcoF_Acof},  $R$ is co-Frobenius if and only if $A=R\#H$ does.
 Now the equivalence of tensor categories given in Lemma 
 \ref{prop: equivalent_cat}  respects dimensions, injective and projective objects, injective hulls and projective covers.
 Thus  assertions $(i)$ to $(x)$ are equivalent by Theorem \ref{thm:characterizations_coFrob}.
\epf

\subsection{Finite coradical filtration}
The converse of Lemma \ref{lem:finfilt_cofrob} holds under stronger hypotheses.

Let $A=R\#H$ as above.
Let $\Pi: A\to R$ the epimorphism of coalgebras given by $\Pi(r\# h)=\varepsilon(h)r$ and $\iota:H\to A$ the canonical injection $h\mapsto1\#h$. 

\begin{lemma}\cite[Thm.~ 11.7.2]{Rad-libro}\label{lem: lemaRad}
If $A_0\iota(H)\subseteq A_0$ and $R_0=\Pi(A_0)$ then $\{\Pi(A_n)\}_{n\geq0}$ is the coradical filtration of $R$.\qed
\end{lemma}

\begin{proposition}\label{prop:cofrob-finite-corad-filt}
	Let $H$ be a co-Fobenius Hopf algebra and
	$R$ a braided co-Frobenius Hopf algebra in $\yd$. If $(R\#H)_0\iota(H)\subseteq (R\#H)_0$ and $R_0=\Pi((R\#H)_0)$ then the coradical filtration of $R$ is finite.
	\end{proposition}
\pf
By Lemma \ref{prop:RcoF_Acof}, $A=R\#H$ is co-Frobenius and therefore by Theorem \ref{thm:characterizations_coFrob} $(xii)$ there exists $n\geq0$ such that $A_n=A$. Thus, by Lemma \ref{lem: lemaRad} $R_n=\Pi(A_n)=\Pi(A)=R$ and the claim follows.
\epf

\end{document}